\documentclass{aptpub}

\authornames{Carl Graham} 
\shorttitle{Chaoticity for multi-class systems and echangeability within classes}

\usepackage{amsfonts,amsmath} 

\newcommand{\EE}{\kern1pt\mathbf{E}\kern1pt}
\newcommand{\PP}{\kern1pt\mathbf{P}\kern1pt}
\newcommand\One{{\mathrm{1} \kern -0.27em \mathrm{I}}}

\renewcommand{\today}{16 October 2008} 

\begin{document}

\title{Chaoticity for multi-class systems\\ and exchangeability within classes}

\author[{\'E}cole Polytechnique, CNRS]{Carl Graham}
\address{CMAP, {\'E}cole Polytechnique, CNRS, 91128 Palaiseau France.}
\email{carl@cmapx.polytechnique.fr}

\begin{abstract} 
Classical results for exchangeable systems of random variables
are extended to multi-class systems satisfying a natural partial exchangeability assumption.
It is proved that
the conditional law of a finite multi-class system,
given the value of the vector of the empirical measures of its classes,
corresponds to \emph{independent} uniform orderings of the samples within \emph{each} class,
and that a family of such systems converges in law \emph{if and only if}
the corresponding
empirical measure vectors converge in law.
As a corollary, 
convergence within \emph{each} class to an infinite i.i.d.\ system implies 
asymptotic independence between \emph{different} classes.
A result implying the Hewitt-Savage 0--1 Law is also extended.
\end{abstract}



\section{Introduction}

Among many others,
Kallenberg~\cite{Kallenberg}, Kingman~\cite{Kingman}, 
Diaconis and Freedman~\cite{Diaconis},
and Aldous~\cite{Aldous} study exchangeable random variables (r.v.)\ with Polish state space.
The related notion of chaoticity (convergence in law to i.i.d.\ random variables) 
appears in many contexts, such as statistical estimation, 
or the asymptotic study of interacting particle systems or communication networks.
It is behind many fruitful heuristics, such as
the ``molecular chaos assumption'' (\emph{Stosszahlansatz}) used by
Ludwig Boltzmann to derive the Boltzmann equation, see
Cercignani \emph{et al.}~\cite[Sect.~2, 4]{Cercignani}.

A sequence of finite exchangeable systems converges
in law to an infinite system if and only if the corresponding sequence 
of empirical measures converges to the 
directing measure of the limit infinite system, given by the de Finetti Theorem.
Hence, chaoticity
is equivalent to the fact that the empirical measures satisfy
a weak law of large numbers, for which A.S. Sznitman developed a compactness-uniqueness method
of proof yielding
propagation of chaos results for varied models of interest.
Sznitman also devised a coupling method for proving chaoticity directly.
See Sznitman~\cite{Sznitman} for a survey, 
and M{\'e}l{\'e}ard~\cite{Meleard:96} and Graham~\cite{Graham:00, Graham:92} for some developments.

The above notions pertain to the study of \emph{similar} random objects, but
many systems in stratified sampling, statistical mechanics, chemistry,
communication networks, biology, etc., involve 
\emph{varied} classes of similar objects (which we call ``particles'').
See for instance Cercignani \emph{et al.}~\cite{Cercignani} 
(``Mixtures'', Subject index p.~454)
and the review papers 
\cite{Bellomo:00, Graham:00, Grunfeld:00, Struckmaier:00}
in a recent book.

Our paper considers natural notions of multi-exchangeability and chaoticity for such multi-class systems,
and extends the above results.
These notions are explicit in Graham~\cite[pp.~78,~81]{Graham:92}, 
and implicit in \cite{Cercignani,Bellomo:00, Grunfeld:00, Struckmaier:00} where
the corresponding limit equations are directly considered.
Graham and Robert~\cite{GrahamRobert} extend Sznitman's coupling method in this context.
For infinite classes, Aldous calls multi-exchangeability ``internal exchangeability''
just before \cite[Corollary~3.9]{Aldous}.

We prove that the conditional law of a finite multi-class system,
given the value of the vector of the empirical measures of its classes,
corresponds to choosing \emph{independent} uniform orderings of the samples within \emph{each} class,
and that a family of such systems converges in law \emph{if and only if}
the corresponding
empirical measure vectors converge in law.
We conclude by extending a result implying the Hewitt-Savage 0--1 Law. 

As a corollary, for a multi-exchangeable system,
chaoticity \emph{within} classes implies asymptotic independence
\emph{between} classes, see Theorem~\ref{main} below.
This striking result allows rigorous derivation
of limit macroscopic models from microscopic dynamics using Sznitman's 
compactness-uniqueness methods, and was a major goal of this paper.

We state as a ``Proposition'' any known result, 
and a ``Theorem'' any result we believe to be new. 
All state spaces $\mathcal{S}$ are Polish,
and the weak topology is used for the space of probability measures $\mathcal{P}(\mathcal{S})$
which is then also Polish, as are products of Polish spaces.
For $k \ge1$ we denote by $\Sigma(k)$ the set of permutations
of $\{1,\ldots,k\}$.

\section{Some classical results}
\label{sclares}

\subsection{Finite and infinite exchangeable systems}

For $N \ge 1$, a finite system $(X^N_n)_{1\le n \le N}$ of random variables (r.v.) 
with state space $\mathcal{S}$ is \emph{exchangeable} if 
\[
\mathcal{L}(X^N_{\sigma(1)}, \ldots, X^N_{\sigma(N)}) 
= \mathcal{L}(X^N_1, \ldots, X^N_N)\,,
\qquad
\forall \sigma \in \Sigma(N)\,.
\]
Then, the conditional law of such a system given the value of its empirical measure
\begin{equation}
  \label{emplaw}
 \Lambda^N = {1\over N} \sum_{n=1}^{N} \delta_{X^N_n} 
\end{equation}
corresponds to a uniform ordering of the $N$ (possibly repeated) values
occurring in $\Lambda^N$
(its atoms, counted according to their multiplicity), see Aldous~\cite[Lemma 5.4 p.~38]{Aldous}.

An \emph{infinite} system $(X_n)_{n \ge 1}$ is \emph{exchangeable} if every
finite subsystem $(X_n)_{1 \le n \le N}$ is exchangeable. 
The de Finetti Theorem, 
see \emph{e.g.} \cite{Kallenberg, Kingman,Diaconis,Aldous},
states that such a system is a mixture
of i.i.d.\ sequences: its law is of the form 
\[
\int P^{\otimes \infty} \mathcal{L}_\Lambda(dP)
\]
where $\mathcal{L}_\Lambda$ is the law of the 
(random) \emph{directing measure} $\Lambda$ 
which can be obtained as
\begin{equation}
\label{dirmeas}
\Lambda = \lim_{N \to \infty} {1\over N} \sum_{n=1}^{N} \delta_{X_n}
\;\;
\textrm{a.s.}
\end{equation}
Thus, laws of infinite exchangeable systems with state space $\mathcal{S}$
and laws of random measures with state space $\mathcal{P}(\mathcal{S})$
are in one-to-one correspondence. 

All this leads to the following fact, see Kallenberg~\cite[Theorem~1.2 p.~24]{Kallenberg} 
and Aldous~\cite[Prop.~7.20 (b) p.~55]{Aldous}.

\begin{proposition}
\label{sznigen}
Let $(X^N_n)_{1\le n \le N}$ for $N \ge1$ be finite exchangeable systems,
and $\Lambda^N$ their empirical measures \eqref{emplaw}.
Then
\[
\lim_{N \to \infty} (X^N_n)_{1\le n \le N} = (X_n)_{n \ge 1} \;\;\textrm{in law},
\]
where the (infinite exchangeable) limit has directing measure $\Lambda$,
if and only if
\[
\lim_{N \to \infty} \Lambda^N = \Lambda
\;\;\textrm{in law}.
\]
\end{proposition}

A sequence $(X^N_n)_{1\le n \le N}$ for $N \ge1$ is $P$-\emph{chaotic}, 
where $P\in \mathcal{P}(\mathcal{S})$, if  
\[
\lim_{N \to \infty} \mathcal{L}(X^N_1, \ldots, X^N_k) = P^{\otimes k}\,,
\qquad
\forall k\ge1\,,
\]
\emph{i.e.}, if it converges in law to an i.i.d.\ system of r.v.\ of law $P$.
The following corollary of Proposition~\ref{sznigen} is proved directly in 
\cite[Prop.~2.2 p.~177]{Sznitman} and \cite[Prop.~4.2 p.~66]{Meleard:96}.

\begin{proposition}
\label{szni}
Let $(X^N_n)_{1\le n \le N}$ for $N \ge1$ be finite exchangeable systems,
$\Lambda^N$ their empirical measures \eqref{emplaw}, and $P\in \mathcal{P}(\mathcal{S})$.
Then, the sequence is $P$-chaotic if and only if
\[
\lim_{N \to \infty} \Lambda^N = P
\;\;\textrm{in law}
\]
and hence in probability, since the limit is deterministic.
\end{proposition}

\subsection{Multi-exchangeable systems}
 
We assume that $C \ge 1$ and state spaces $\mathcal{S}_i$ for $1 \le i \le C$ are fixed.
For a multi-index $\mathbf{N}=(N_i)_{1 \le i \le C} \in \mathbb{N}^C$ 
we consider a multi-class system
\begin{equation}
\label{eq:mcys}
(X^\mathbf{N}_{n,i})_{1\le n \le N_i,\, 1\le i \le C}\,,
\qquad
X^\mathbf{N}_{n,i}
\textrm{ with state space $\mathcal{S}_i$},
\end{equation}
where $X^\mathbf{N}_{n,i}$ is the $n$-th particle, or object, of class $i$, and
say that it is \emph{multi-exchangeable} if its law is invariant under
permutation of the particles \emph{within} classes:
\[
\mathcal{L}\bigl( (X^\mathbf{N}_{\sigma_i(n),i})_{1\le n \le N_i,\, 1\le i \le C} \bigr)
=
\mathcal{L}\bigl( (X^\mathbf{N}_{n,i})_{1\le n \le N_i,\, 1\le i \le C} \bigr)\,,
\qquad
\forall \sigma_i \in \Sigma(N_i)\,.
\]
This natural assumption means that particles of a class are
statistically indistinguishable, and obviously
implies that $(X^\mathbf{N}_{n,i})_{1\le n \le N_i}$ is exchangeable for $1\le i \le C$.
It is sufficient to check that it is true 
when all $\sigma_i$ but one are the identity.
We introduce the \emph{empirical measure vector}, with samples in 
$\mathcal{P}(\mathcal{S}_1) \times \cdots \times \mathcal{P}(\mathcal{S}_C)$,
\begin{equation}
  \label{eq:emmeve}
(\Lambda^\mathbf{N}_i)_{1\le i \le C}\,,
\qquad
\Lambda^\mathbf{N}_i = {1 \over N_i}\sum_{n=1}^{N_i} \delta_{X^\mathbf{N}_{n,i}}\,.
\end{equation}

We say that the multi-class system $(X_{n,i})_{n \ge 1, 1\le i \le C}$
with \emph{infinite} classes
is \emph{multi-exchangeable} if every finite sub-system $(X_{n,i})_{1 \le n \le N_i,\, 1\le i \le C}$
is multi-exchangeable.
Particles of class~$i$ form an exchangeable system, which has a directing measure $\Lambda_i$, and we call
$(\Lambda_i)_{1 \le i \le C}$ the \emph{directing measure vector}.

The following result is given in Aldous~\cite[Cor.~3.9 p.~25]{Aldous}
and  attributed to
de Finetti.
A remarkable fact is conditional independence between \emph{different} classes.

\begin{proposition}
\label{conind}
Let $(X_{n,i})_{n \ge 1,\, 1\le i \le C}$ be an infinite multi-exchangeable system,
and $\Lambda_i$ be the directing measure of $(X_{n,i})_{n \ge 1}$.
Given the directing measure vector $(\Lambda_i)_{1 \le i \le C}$, the 
$X_{n,i}$ for $n \ge 1$ and $1\le i \le C$ are conditionally independent, 
and
$X_{n,i}$ has conditional law $\Lambda_i$.
\end{proposition}

\section{The extended results}

We shall extend to multi-exchangeable systems the main results for exchangeable systems, 
which hold even though the symmetry assumption and resulting structure is much weaker.
Indeed,
the symmetry order of the multi-exchangeable system \eqref{eq:mcys}
is $N_1! \cdots N_C!$ whereas the symmetry order of an exchangeable system of same size
is the much larger
$(N_1 + \cdots + N_C)!$.

The following extension of \cite[Lemma~5.4 p.~38]{Aldous} 
(stated in words at the beginning of Section~\ref{sclares}) 
shows that, for a finite multi-exchangeable system,
the classes are \emph{conditionally independent}
given the vector of the empirical measures within each class.
Hence, \emph{no further information}
can be attained on its law
by cleverly trying to involve what happens for different classes.

A statistical interpretation of this remarkable fact is that the empirical measure vector is a 
\emph{sufficient statistic} for the law of the system,
the family of all such laws being trivially parameterized by the laws themselves.

\begin{theorem}
\label{suffstat}
Let $(X^\mathbf{N}_{n,i})_{1\le n \le N_i,\, 1\le i \le C}$ be a finite multi-exchangeable system
as in \eqref{eq:mcys}. 
Then its conditional law, given the value of the empirical measure vector 
$(\Lambda^\mathbf{N}_i)_{1\le i \le C}$ defined in \eqref{eq:emmeve},
corresponds to \emph{independent} uniform orderings for $1 \le i \le C$ 
of the $N_i$ values of the particles of class~$i$ (possibly repeated), 
which are the atoms of the value of $\Lambda^\mathbf{N}_i$
(counted with their multiplicities).
\end{theorem}

\begin{proof}
Multi-exchangeability and the obvious fact that
\[
\Lambda^\mathbf{N}_j 
= {1 \over N_j}\sum_{n=1}^{N_j} \delta_{X^\mathbf{N}_{n,j}}
= {1 \over N_j}\sum_{n=1}^{N_j} \delta_{X^\mathbf{N}_{\sigma(n),j}}\,,
\qquad
\forall \sigma \in \Sigma(N_j)\,,
\]
imply that for all
$g : \mathcal{P}(\mathcal{S}_1) \times \cdots \times \mathcal{P}(\mathcal{S}_C) \rightarrow \mathbb{R}_+$ 
and $f_i : \mathcal{S}_i^{N_i} \rightarrow \mathbb{R}_+$
we have
\begin{eqnarray}
\label{cenfor}
&&\kern-6.5mm
\EE\left[
g\bigl((\Lambda^\mathbf{N}_j)_{1 \le j \le C}\bigr)
\prod_{i=1}^C 
f_i(X^\mathbf{N}_{1,i}, \ldots X^\mathbf{N}_{N_i,i})
\right]
\nonumber\\
&&\kern-5.5mm{}
= {1 \over N_1 !}\sum_{\sigma_1 \in \Sigma(N_1)} \cdots {1 \over N_C !}\sum_{\sigma_C \in \Sigma(N_C)}
\EE\left[
g\bigl((\Lambda^\mathbf{N}_j)_{1 \le j \le C}\bigr)
\prod_{i=1}^C 
f_i(X^\mathbf{N}_{\sigma_i(1),i}, \ldots X^\mathbf{N}_{\sigma_i(N_i),i})
\right]
\nonumber\\
&&\kern-5.5mm{}
=
\EE\!\left[
g\bigl((\Lambda^\mathbf{N}_j)_{1 \le j \le C}\bigr)
\prod_{i=1}^C  
{1 \over N_i !} 
\sum_{\sigma \in \Sigma(N_i) }
f_i(X^\mathbf{N}_{\sigma(1),i}, \ldots, X^\mathbf{N}_{\sigma(N_i),i})
\right]
\nonumber\\
&&\kern-5.5mm{}
=
\EE\!\left[
g\bigl((\Lambda^\mathbf{N}_j)_{1 \le j \le C}\bigr)
\prod_{i=1}^C  
\left\langle f_i,
{1 \over N_i !} 
\sum_{\sigma \in \Sigma(N_i) }
\delta_{ (X^\mathbf{N}_{\sigma(1),i}, \ldots, X^\mathbf{N}_{\sigma(N_i),i}) }
\right\rangle
\right]
\end{eqnarray}
where the empirical measure 
\[
{1 \over N_i !} 
\sum_{\sigma \in \Sigma(N_i) }
\delta_{ (X^\mathbf{N}_{\sigma(1),i}, \ldots, X^\mathbf{N}_{\sigma(N_i),i}) }
\]
corresponds to exhaustive uniform draws without replacement among 
the atoms 
$X^\mathbf{N}_{1,i}$, \dots\,, $X^\mathbf{N}_{N_i,i}$ of $\Lambda^\mathbf{N}_i$
counted according to multiplicity, and hence is a function of $\Lambda^\mathbf{N}_i$.
Since $g$ is arbitrary, the characteristic property of conditional expectation
yields that
\[
\EE\Biggl[
\prod_{i=1}^C 
f_i(X^\mathbf{N}_{1,i}, \ldots X^\mathbf{N}_{N_i,i})
\,\bigg |\,(\Lambda^\mathbf{N}_i)_{1 \le i \le C}
\Biggr]
= \prod_{i=1}^C  
\left\langle f_i,
{1 \over N_i !} 
\sum_{\sigma \in \Sigma(N_i) }
\delta_{ (X^\mathbf{N}_{\sigma(1),i}, \ldots, X^\mathbf{N}_{\sigma(N_i),i}) }
\right\rangle
\]
which finishes the proof, since the $f_i$ are arbitrary and the spaces Polish.
\end{proof}

This result and Proposition~\ref{conind} lead to the following extension
of Proposition~\ref{sznigen}. 
We denote by $\lim_{\mathbf{N} \to \infty}$ the limit along a fixed arbitrary 
subsequence of $\mathbf{N} \in \mathbb{N}^C$ such that $\min_{1 \le i \le C} N_i$ goes to infinity.

\begin{theorem}
\label{maingen}
We consider a family of finite  multi-exchangeable multi-class systems
\[
(X^\mathbf{N}_{n,i})_{1\le n \le N_i,\, 1\le i \le C}\,,
\qquad
\mathbf{N} \in \mathbb{N}^C\,, 
\]
all of the form
\eqref{eq:mcys} with the same $C\ge1$ and state spaces $\mathcal{S}_i$,
and the corresponding empirical measure vectors $(\Lambda^\mathbf{N}_i)_{1\le i \le C}$
given in \eqref{eq:emmeve}.
Then
\[
\lim_{\mathbf{N} \to \infty} (X^\mathbf{N}_{n,i})_{1\le n \le N_i,\, 1\le i \le C}
= (X_{n,i})_{n \ge 1,\, 1\le i \le C}
\;\; \textrm{in law},  
\]
where the (infinite multi-exchangeable) limit 
has directing measure vector $(\Lambda_i)_{1 \le i \le C}$,
if and only if
\[
\lim_{\mathbf{N} \to \infty} (\Lambda^\mathbf{N}_i)_{1 \le i \le C} = (\Lambda_i)_{1 \le i \le C}
\;\; \textrm{in law.}
\]
\end{theorem}

\begin{proof}
Since the state spaces are Polish, it is enough to prove that
for arbitrary $k \ge1$ and bounded continuous $f_i : \mathcal{S}_i^k \rightarrow \mathbb{R}$ 
for $1 \le i \le C$
we have
\begin{equation}
\label{oneway}
\lim_{\mathbf{N} \to \infty} 
\EE\left[
\prod_{i=1}^C 
f_i(X^\mathbf{N}_{1,i}, \ldots X^\mathbf{N}_{k,i})
\right]
= 
\EE\left[
\prod_{i=1}^C 
f_i(X_{1,i}, \ldots X_{k,i})
\right]
\end{equation}
if and only if
\begin{equation}
\label{theother}
\lim_{\mathbf{N} \to \infty} 
\EE\!\left[
\prod_{i=1}^C  
\left\langle f_i, (\Lambda^\mathbf{N}_i)^{\otimes k} \right\rangle
\right]
= 
\EE\!\left[
\prod_{i=1}^C  
\left\langle f_i, \Lambda_i^{\otimes k} \right\rangle
\right].
\end{equation}
Let $(m)_k = \frac{m!}{(m-k)!} = m(m-1)\cdots(m-k+1)$ for $m\ge 1$
and, for $N_i \ge k$,
\[
\Lambda^{\mathbf{N},k}_i 
=
{1 \over (N_i)_k} 
\sum_{ \substack{ 1 \le n_1,\ldots, n_k \le N_i \\ \textrm{distinct} }}
\delta_{ (X^\mathbf{N}_{n_1,i}, \ldots, X^\mathbf{N}_{n_k,i}) } 
\]
denote the empirical measure for distinct $k$-tuples in class~$i$,
corresponding to sampling $k$ times \emph{without} 
replacement among $X^\mathbf{N}_{1,i}$, \dots\,, $X^\mathbf{N}_{N_i,i}$.
Theorem~\ref{suffstat} implies that
\begin{eqnarray}
\label{prf1}
\EE\left[
\prod_{i=1}^C 
f_i(X^\mathbf{N}_{1,i}, \ldots X^\mathbf{N}_{k,i})
\right]
&=&
\EE\left[
\EE\left[
\prod_{i=1}^C 
f_i(X^\mathbf{N}_{1,i}, \ldots X^\mathbf{N}_{k,i})
\,\bigg|\, (\Lambda^\mathbf{N}_i)_{1 \le i \le C}
\right]
\right]
\nonumber \\
&=&
\EE\!\left[
\prod_{i=1}^C  
\left\langle f_i,
\Lambda^{\mathbf{N},k}_i 
\right\rangle
\right]
\quad
\end{eqnarray}
(which follows directly 
from \eqref{cenfor} with $g=1$ and the extensions of $f_i$ on $\mathcal{S}_i^{N_i}$) 
and Proposition~\ref{conind} similarly implies that
\begin{equation}
\label{conddec}
\EE\left[
\prod_{i=1}^C 
f_i(X_{1,i}, \ldots X_{k,i})
\right]
=
\EE\!\left[
\prod_{i=1}^C  
\left\langle f_i, \Lambda_i^{\otimes k} \right\rangle
\right]\,.
\end{equation}
The corresponding empirical measure for sampling \emph{with} replacement is given by
\begin{eqnarray*}
(\Lambda^\mathbf{N}_i)^{\otimes k}
&=&
{1 \over N_i^k} \sum_{1 \le n_1,\ldots, n_k \le N_i}
\delta_{ (X^\mathbf{N}_{n_1,i}, \ldots, X^\mathbf{N}_{n_k,i}) }
\\
&=& 
\frac{(N_i)_k}{N_i^k}\,
\Lambda^{\mathbf{N},k}_i 
+ 
{1 \over N_i^k} \sum_{ \substack{ 1 \le n_1,\ldots, n_k \le N_i \\ \textrm{not distinct} }}
\delta_{ (X^\mathbf{N}_{n_1,i}, \ldots, X^\mathbf{N}_{n_k,i}) }
\end{eqnarray*}
and in total variation norm
$\Vert \mu \Vert = \sup\{\,\langle \phi, \mu \rangle : \Vert \phi \Vert_\infty \le 1\,\}$
we have
\begin{equation}
\label{tveq}
\left\Vert 
(\Lambda^\mathbf{N}_i)^{\otimes k}
-
\Lambda^{\mathbf{N},k}_i 
\right\Vert
\le 2 {N_i^k - (N_i)_k  \over N_i^k}  \le {k(k-1) \over N_i}
\end{equation}
where we bound $N_i^k - (N_i)_k$ by counting $k(k-1)/2$ possible positions for two
identical indices with $N_i$ choices and $N_i^{k-2}$ choices for the other $k-2$ positions.
Hence,
if
\eqref{oneway}
holds then \eqref{prf1}, \eqref{tveq} and \eqref{conddec} imply \eqref{theother},
and conversely,
if \eqref{theother} holds then \eqref{tveq}, \eqref{prf1} and \eqref{conddec} imply
\eqref{oneway},
which concludes the proof.
\end{proof}

Let $P_i \in \mathcal{P}(\mathcal{S}_i)$ for $1\le i \le C$.
We say that the family of finite multi-class systems
such as in Theorem~\ref{maingen}
is $(P_1,\ldots, P_C)$-\emph{chaotic} if 
\[
\lim_{\mathbf{N} \to \infty}
\mathcal{L}\bigl( (X^\mathbf{N}_{n,i})_{1\le n \le k,\, 1\le i \le C} \bigr)
= P_1^{\otimes k}\otimes \cdots \otimes P_C^{\otimes k}\,,
\qquad
 \forall k\ge1\,.
\]
This means that the multi-class systems converge to an \emph{independent}
system, in which particles of class~$i$ have law $P_i$. 
We state a striking corollary of Theorem~\ref{maingen}. 
 
\begin{theorem}
\label{main}
We consider a family  of finite multi-exchangeable 
multi-class systems such as in Theorem~\ref{maingen},
and $P_i \in \mathcal{P}(\mathcal{S}_i)$ for $1\le i \le C$.
Then the family is $(P_1,\ldots, P_C)$-chaotic
if and only if the
$(X^\mathbf{N}_{n,i})_{1\le n \le N_i}$ are $P_i$-chaotic for $1\le i \le C$.
\end{theorem}

\begin{proof}
Since $(P_i)_{1 \le i \le C}$ is deterministic, 
$\lim_{\mathbf{N} \to \infty}(\Lambda^\mathbf{N}_i)_{1 \le i \le C} = (P_i)_{1 \le i \le C}$
in law if and only if 
$\lim_{\mathbf{N} \to \infty}\Lambda^\mathbf{N}_i = P_i$ 
in law for $1\le i \le C$. 
We conclude using Theorem~\ref{maingen}.
\end{proof}

We finish with the following extension of Aldous~\cite[Cor.~3.10 p.~26]{Aldous}
and of the Hewitt-Savage 0--1 Law.
For $k \ge 1$, we say that a set 
\[
B \subset \mathcal{S}_1^{\infty} \times \cdots \times \mathcal{S}_C^{\infty}
\]
is $k$-multi-exchangeable if 
for all permutations $\sigma_i$ 
of $\{1,2,\ldots\}$ leaving $\{k+1,k+2,\ldots\}$ invariant, $1 \le i \le C$, we have
\[
(x_{n,i})_{n \ge 1,\, 1\le i \le C} \in B 
\Leftrightarrow
(x_{\sigma_i(n),i})_{n \ge 1,\, 1\le i \le C}  \in B\,.
\]
We define the multi-exchangeable $\sigma$-algebra
\[
\mathcal{E} = \bigcap_{k \ge 1} \mathcal{E}_k\,,
\qquad
\mathcal{E}_k = \bigl\{
\{ (X_{n,i})_{n \ge 1 ,\, 1\le i \le C} \in B \}
: B \textrm{ is $k$-multi-exchangeable} \bigr\} \,,
\]
and multi-tail $\sigma$-algebra 
\[
\mathcal{T} = \bigcap_{k \ge 1} \mathcal{T}_k\,,
\qquad
\mathcal{T}_k = \sigma\bigl((X_{n,i})_{n \ge k ,\, 1\le i \le C}\bigr)\,.
\]
Clearly, 
$\mathcal{T}_{k+1} \subset \mathcal{E}_k$ and hence $\mathcal{T} \subset \mathcal{E}$.

\begin{theorem}
Let $(X_{n,i})_{n \ge 1,\, 1\le i \le C}$ be an infinite multi-exchangeable system
with directing measure vector $(\Lambda_i)_{1 \le i \le C}$. Then
\[
\sigma((\Lambda_i)_{1 \le i \le C}) = \mathcal{T} = \mathcal{E}
\;\;\textrm{a.s.}
\]
If moreover the $X_{n,i}$ are independent, then $P(A) \in \{0,1\}$
for all $A \in \mathcal{E}$.
\end{theorem}

\begin{proof}
Consideration of \eqref{dirmeas} yields 
$\sigma((\Lambda_i)_{1 \le i \le C}) \subset \mathcal{T}$,
a.s., and we have seen that $\mathcal{T} \subset \mathcal{E}$, 
hence the first statement is true if
$\mathcal{E} \subset \sigma((\Lambda_i)_{1 \le i \le C})$, a.s.
Now, let $A \in \mathcal{E}$. 
For every $k\ge1$, since $A \in \mathcal{E}_k$, 
there is some $k$-multi-exchangeable set $B_k$ such that
\[
A = \{ (X_{n,j})_{n \ge 1 ,\, 1\le j \le C} \in B_k \}
\]
and hence,
for all permutations $\sigma_i$ 
of $\{1,2,\ldots\}$ leaving $\{k+1,k+2,\ldots\}$ invariant, 
\begin{eqnarray*}
(\mathbf{1}_{A}, X_{\sigma_i(n),i})_{n \ge 1,\, 1\le i \le C}
&=& 
(\mathbf{1}_{B_k}((X_{n,j})_{n \ge 1 ,\, 1\le j \le C}), X_{\sigma_i(n),i}
)_{n \ge 1,\, 1\le i \le C}  
\\
&=& 
(
\mathbf{1}_{B_k}((X_{\sigma_j(n),j})_{n \ge 1 ,\, 1\le j \le C}), X_{\sigma_i(n),i}
)_{n \ge 1,\, 1\le i \le C}  
\end{eqnarray*}
and the multi-exchangeability of $(X_{n,i})_{n \ge 1,\, 1\le i \le C}$ implies that
\begin{eqnarray*}
\mathcal{L}\bigl((\mathbf{1}_A, X_{\sigma_i(n),i})_{n \ge 1,\, 1\le i \le C}\bigr)
&=& 
\mathcal{L}\bigl(
(\mathbf{1}_{B_k}((X_{\sigma_j(n),j})_{n \ge 1 ,\, 1\le j \le C}), X_{\sigma_i(n),i}
\bigr)_{n \ge 1,\, 1\le i \le C}
\bigr)
\\
&=& 
\mathcal{L}\bigl(
(\mathbf{1}_{B_k}((X_{n,j})_{n \ge 1 ,\, 1\le j \le C}), X_{n,i}
\bigr)_{n \ge 1,\, 1\le i \le C}
\bigr)
\\
&=& 
\mathcal{L}\bigl((\mathbf{1}_A, X_{n,i})_{n \ge 1,\, 1\le i \le C}\bigr)\,.
\end{eqnarray*}
Thus
$(\mathbf{1}_A, X_{n,i})_{n \ge 1,\, 1\le i \le C}$ is infinite multi-exchangeable, and
Proposition~\ref{conind} implies that
the $(\mathbf{1}_A, X_{n,i})$ are conditionally independent given $(\hat\Lambda_i)_{1 \le i \le C}$
and have conditional laws $\hat\Lambda_i$, 
where considering \eqref{dirmeas} we have
\[
\hat\Lambda_i = \lim_{N \to \infty} {1 \over N} \sum_{n=1}^N \delta_{ (\mathbf{1}_A, X_{n,i}) }
= \delta_{\mathbf{1}_A} \otimes 
\lim_{N \to \infty} {1 \over N} \sum_{n=1}^N \delta_{X_{n,i}}
= \delta_{\mathbf{1}_A} \otimes \Lambda_i
\;\;\textrm{a.s.}
\]
Hence, for arbitrary $k \ge1$ and Borel sets $B_{n,i} \subset \mathcal{S}_i$ 
for $1 \le n \le k$ and $1 \le i \le C$,
\[
\displaylines{\quad
\PP\bigl(X_{n,i} \in B_{n,i} : 1 \le n \le k, 1 \le i \le C 
\,\big|\, A, (\Lambda_i)_{1 \le i \le C}\bigr)
\hfill\cr\hfill
= \PP\bigl(X_{n,i} \in B_{n,i} : 1 \le n \le k, 1 \le i \le C 
\,\big|\, (\hat\Lambda_i)_{1 \le i \le C}\bigr)
= \prod_{1 \le n \le k, 1 \le i \le C} \Lambda_i(B_{n,i})
\quad}
\]
is a function of $(\Lambda_i)_{1 \le i \le C}$, 
conditionally to which
$A$ and $(X_{n,i})_{n \ge 1,\, 1\le i \le C}$
are thus independent.
Since $A \in \mathcal{E}$ is arbitrary, we deduce that
$\mathcal{E} \subset \sigma((X_{n,i})_{n \ge 1,\, 1\le i \le C})$
and $(X_{n,i})_{n \ge 1,\, 1\le i \le C}$
are conditionally independent given $(\Lambda_i)_{1 \le i \le C}$, 
which implies $\mathcal{E} \subset \sigma((\Lambda_i)_{1 \le i \le C})$ a.s.
This proves the first statement, from which the second follows 
since $\mathcal{T}$ is a.s.\ trivial if the
$X_{n,i}$ are independent, see the Kolmogorov 0--1 Law.
\end{proof}

\section{Concluding remarks}

The important bound \eqref{tveq} is a combinatorial estimate of the difference between 
sampling with and without replacement,
see Aldous~\cite[Prop.~5.6 p.~39]{Aldous} and Diaconis and Freedman~\cite[Theorem~13 p.~749]{Diaconis} 
for related results.
It is used in \cite{Diaconis} to prove the de Finetti Theorem.

Theorem~\ref{main} allows proving $(P_1,\ldots, P_C)$-chaoticity results
by use of Proposition~\ref{szni} and Sznitman's compactness-uniqueness methods
for proof that the empirical measures $\Lambda^\mathbf{N}_i$
converge in law to $P_i$ for $1 \le i \le C$. 
This was the main motivation for this paper,
as can be seen by its title. In the reviewing process, the referee's suggestions lead to 
a much improved and fuller study of multi-exchangeable systems.

The techniques developed in this paper could also extend convergence results, such as
Kallenberg~\cite[Theorem~1.3 p.~25]{Kallenberg} and Aldous~\cite[Prop.~7.20 (a) p.~55]{Aldous},
suited for a family of multi-exchangeable systems of fixed possibly infinite class sizes
depending on a parameter.
We refrain do to so for the sake of coherence.

\section*{Acknowledgments}

We thank the anonymous referee for a careful reading 
and many stimulating remarks. He helped us rediscover the depth and beauty 
of Aldous's treatise~\cite{Aldous}.


\end{document}